\documentclass[12pt,letterpaper]{amsart}
\usepackage{geometry}
\usepackage{mathrsfs}
\usepackage{amsfonts}
\usepackage{amsthm}
\usepackage{amssymb,latexsym,amscd} 
\usepackage[dvips]{color}

\usepackage{hyperref}
\newtheorem{theorem}{Theorem}

\newtheorem{corollary}{Corollary}
\theoremstyle{remark}
\newtheorem{remark}{Remark}

\newtheorem{proposition}{Proposition}

\newcommand{\C}{\mathbb{C}}

\newcommand{\D}{\Omega}
\newcommand{\dbar}{\overline{\partial}}

\title[Irregularity of the $\dbar$-Neumann problem on non-smooth domains]{A remark On Irregularity of the $\dbar$-Neumann problem on non-smooth domains}
\author{S\"onmez \c Sahuto\u glu}
\address{Department of Mathematics \\ University of Michigan \\
Ann Arbor, MI 48109-1043}

\email{sonmez@umich.edu}
\date{\today}
\begin{document}

\begin{abstract}
It is an observation due to J.J. Kohn that  for a smooth bounded pseudoconvex domain $\D$ in $\C^n$ there exists $s>0$ such that the $\dbar$-Neumann operator on $\D$  maps $W^s_{(0,1)}(\D)$ (the space of $(0,1)$-forms with coefficient functions in $L^2$-Sobolev space of order $s$) into itself continuously. We show that this conclusion does not hold  without the smoothness assumption by constructing a bounded pseudoconvex
domain $\Omega$ in $\mathbb{C}^{2}$, smooth except at one point, whose $\dbar$-Neumann operator is  not bounded on $W^s_{(0,1)}(\D)$ for any $s>0$.
\end{abstract}

\maketitle

Let $W^s(\D)$ and $W^s_{(p,q)}(\D)$ denote the $L^2$-Sobolev space on $\D$ of order $s$ and the space of $(p,q)$-forms with coefficient functions in $W^s(\D)$, respectively. Also $\|.\|_{s,\D}$ denotes the norms on $W^s_{(p,q)}(\D)$. Let $N_q$ denote the inverse of the complex Laplacian, $\dbar\dbar^*+\dbar^*\dbar$, on square integrable $(0,q)$-forms. It is an observation of Kohn, as the following proposition says,  that on a smooth bounded pseudoconvex domain the $\dbar$-Neumann problem is regular in the Sobolev scale for sufficiently small levels.  

\begin{proposition}[Kohn]
Let $\D$ be a smooth bounded pseudoconvex domain in $\C^n$. There exist positive $\varepsilon$ and $C$ (depending on $\D$) such that 
$$\|N_q u\|_{\varepsilon,\D}\leq C\|u\|_{\varepsilon,\D},\|\dbar N_q u\|_{\varepsilon,\D}\leq C\|u\|_{\varepsilon,\D}, \|\dbar^* N_q u\|_{\varepsilon.\D}\leq C\|u\|_{\varepsilon,\D}$$ for $u\in W^s_{(0,q)}(\D)$ and $1\leq  q \leq n$.
\end{proposition}
 We show that if one drops the smoothness assumption then the  $\dbar$-Neumann operator, $N_1,$ may not map any positive Sobolev space into itself continuously. 

\begin{theorem}\label{thm1}
There exists a bounded pseudoconvex domain $\D$ in $\C^2$, smooth except one point, such that the $\dbar$-Neumann operator on $\D$ is not bounded on $W^s_{(0,1)}(\D)$ for any $s>0$.
\end{theorem}

\begin{proof} We will build the domain by attaching infinitely
many worm domains (constructed by Diederich and Forn\ae{}ss in \cite{DF}) 
with progressively larger winding.
Let $\D_j$ be a worm domain, a smooth bounded pseudoconvex domain, in $\C^2$ that winds $2\pi j$  such that 
$$\D_j\subset \{(z,w)\in \C^2:|z|<2^{-j},4^{-j}<|w|<4^{-j}2\}$$ 
for $j=1,2,\ldots$. Let $\gamma_j$ be a straight line that connects an extreme point on the cap of $\D_j$ to a closest point on the cap of $\D_{j+1}$. Then using the barbell lemma (see \cite{FS,HW}) we get a bounded pseudoconvex domain $\D$ that is smooth except one point $(0,0)\in b\D$. Notice that $\D$ is the union of $\D_j\subset \D$ for $j=1,2,\ldots$ and all connecting bands. In the rest of the proof we will show that if the $\dbar$-Neumann operator on $\D$ is continuous on $W^s_{(0,1)}(\D)$ then the $\dbar$-Neumann operator on $\D_j$ is continuous on  $W^s_{(0,1)}(\D_j)$ for $j=1,2,\ldots$. However this is a contradiction with  a theorem of Barrett(\cite{B92}).  Let us define $\Box^j=\dbar\dbar^*+\dbar^*\dbar$ on $L^2_{(0,1)}(\D_j)$, and  $\Box=\dbar\dbar^*+\dbar^*\dbar$ on $L^2_{(0,1)}(\D)$.  Let us fix $j$ and choose a defining function 
$\rho$ for $\D_{j}$ such that $\|\nabla \rho\|=1$ on $b\D_{j}$. Let $\nu=Re\left(\sum_{j=1}^{2}\frac{\partial \rho}{\partial \bar z_{j}}\frac{\partial }{\partial z_{j}}\right)$ and $J$ denote  the complex structure of $\C^2$. Now we will construct a smooth cut off function that fixes the domain of $\Box$ and $\Box^j$ under multiplication. We can choose open sets $U_1,U_2,$ and $U_3$ and $\chi\in C_0^{\infty}( U_2)$ such that 
\begin{itemize}
\item[i)] $U_1\subset \subset U_2 \subset\subset U_3$,
\item[ii)] $U_{1},U_{2},$ and $U_3$ contain all boundary points of $\D_j$ that meet the (strongly pseudoconvex) band created using $\gamma_j$ and $\gamma_{j-1}$, and they do not contain any weakly pseudoconvex boundary point of $\D_j$,
\item[iii)] $0\leq \chi \leq 1$, $\chi \equiv 1$ on $U_1$,
\item[iv)] there exists an open set $U$ such that $b\D_j\cup U_2\subset\subset U$ and the following two ordinary differential equations can be solved in $U$ 
\begin{eqnarray}
&\nu(\widetilde{\psi})=0,& \widetilde{\psi}|_{b\D_j}=\chi,\\
&\nu(\widetilde{\phi})=-J(\nu)(\chi),& \widetilde{\phi}|_{b\D_j}=0.
\end{eqnarray}
\end{itemize} 
Notice that $\widetilde{\psi} \equiv 1$ and $\widetilde{\phi}\equiv 0$ on $U_1$, and $\widetilde{\psi}= \widetilde{\phi}= 0$ in a neighborhood of the set of weakly pseudoconvex boundary points of $\D_j$. We choose a neighborhood $V\subset\subset U$ of $b\D_j$ and $\widetilde{\chi}\in C^{\infty}_0(V)$ such that $\widetilde{\chi}\equiv 1 $ in a neighborhood $\widetilde{V}$ of $b\D_j$. Let us define  $\phi=\widetilde{\chi}\widetilde{\phi},\psi=\widetilde{\chi}\widetilde{\psi},$ and $\xi=\psi+i\phi$. We like to make some observation about $\xi$ that will be useful later:
\begin{itemize}
\item[i)] $\xi \equiv 1$ on $\widetilde{V}\cap U_1$,
\item[ii)] $(\nu+iJ(\nu))(\xi)\equiv 0$ on $b\D_j$,
\item[iii)] $\xi\equiv 0$ in a neighborhood of the weakly  pseudoconvex boundary points of $\D_j$.
\end{itemize}

\smallskip
\noindent{\it{Claim: If $f\in Dom(\Box^j)$ then $\xi f\in Dom(\Box^j)$ and $(1-\xi)f\in Dom(\Box)$.}}

\smallskip
\noindent{\it{Proof of Claim: }} First we will show that $\xi f\in Dom(\Box^j)$ then we will talk about  how one can show that $(1-\xi)f\in Dom(\Box)$. 

One can easily show that $\xi f\in Dom(\dbar^*)\cap Dom(\dbar)$ (on $\D_j$). On the other hand, by Kohn-Morrey-H\"ormander formula \cite{shaw} since the $L^2$-norms of  any ``bar" derivatives of any terms of $f$ on $\D_j$ is dominated by $\|\dbar f\|_{\D_j}+\| \dbar^* f\|_{\D_j}$  we have $\dbar^*(\xi f)\in Dom(\dbar)$. So we need to show that $\dbar (\xi f)=\dbar\xi\wedge f+\xi\dbar f\in Dom(\dbar^*)$. Since $\xi \dbar f\in Dom(\dbar^*)$ we only need to show that $\dbar\xi\wedge f \in Dom(\dbar^*)$. We will use the special boundary frames. Let $$L_{\tau}=\frac{\partial \rho}{\partial z_{1}}\frac{\partial }{\partial z_{2}}-\frac{\partial \rho}{\partial z_{2}}\frac{\partial }{\partial z_{1}},\quad  L_{\nu}=\frac{\partial \rho}{\partial \bar z_{1}}\frac{\partial }{\partial z_{1}}+\frac{\partial \rho}{\partial \bar z_{2}}\frac{\partial }{\partial z_{2}}.$$  Also let $w_{\tau}$ and $w_{\nu}$ be the dual  $(1,0)$-forms. We note that $L_{\nu}=\nu-iJ(\nu)$ and so $\overline{L}_{\nu}(\xi)\equiv 0$ on $b\D_j$.  We can write $f=f_{\tau}\overline{w}_{\tau}+f_{\nu}\overline{w}_{\nu}$. Therefore, $\dbar\xi\wedge f = (\overline{L}_{\tau}(\xi)f_{\nu}-\overline{L}_{\nu}(\xi)f_{\tau})\overline{w}_{\tau}\wedge \overline{w}_{\nu}$.  Using the fact that  $f_{\nu}\in W_0^1(\D_j)$ (it is easy to see this for $f\in C^1(\overline{\D}_j)$. For $f\in  Dom(\dbar^*)\cap Dom(\dbar)$ one can use the fact that $\Delta:W^1_0(\D_j)\to W^{-1}(\D_j)$ is an isomorphism and the density lemma \cite[Lemma 4.3.2]{shaw} to see this) and $\overline{L}_{\tau}(\xi)$ is smooth we may reduce the problem of showing $\dbar\xi\wedge f\in Dom(\dbar^*)$ to show the following  
$$\overline{L}_{\nu}(\xi)f_{\tau}\overline{w}_{\tau}\wedge \overline{w}_{\nu} \in Dom(\dbar^*).$$ 
Let $\{\phi_k\}_{k=1}^{\infty}$ be a sequence of smooth compactly supported functions converging to $\overline{L}_{\nu}(\xi)$ in $C^1$-norm and $u$ be a $(0,1)$-form with smooth compactly supported coefficient functions in $\D_j$. Then  
$$\langle\overline{L}_{\nu}(\xi)f_{\tau}\overline{w}_{\tau}\wedge \overline{w}_{\nu},\dbar u\rangle_{\D_j}=\lim_{k\to \infty}\langle \phi_kf_{\tau}\overline{w}_{\tau}\wedge \overline{w}_{\nu},\dbar u \rangle_{\D_j} $$ where $\langle,\rangle_{\D_j}$ is the inner product on forms on $\D_j$. If we integrate by parts and use  $\lim_{k\to\infty} \| L_l(\phi_kf_{\tau})\|_{\D_j}=\| L_l(\overline{L}_{\nu}(\xi )f_{\tau})\|_{\D_j}$ for $l=\tau,\nu$ we can reduce the problem of showing $\dbar\xi\wedge f\in Dom(\dbar^*)$ to showing that $\| \frac{\partial}{\partial z_1}(\overline{L}_{\nu}(\xi )f_{\tau})\|_{\D_j}$ and $\| \frac{\partial}{\partial z_2}(\overline{L}_{\nu}(\xi )f_{\tau})\|_{\D_j}$  are finite. One can show that
\begin{equation*}
\left\|\frac{\partial}{\partial z_m}(\overline{L}_{\nu}(\xi)f_{\tau}) \right\|_{\D_j}=\lim_{k\to \infty}\left\|\frac{\partial}{\partial z_m}(\phi_kf_{\tau}) \right\|_{\D_j}
=\lim_{k\to \infty}\left\|\frac{\partial}{\partial \bar z_m}(\phi_kf_{\tau}) \right\|_{\D_j}.
\end{equation*}
On the second equality we used integration by parts. On the other hand,  we have
\begin{eqnarray*}
\lim_{k\to \infty}\left\|\frac{\partial}{\partial \bar z_m}(\phi_kf_{\tau}) \right\|_{\D_j}
&=& \left\|\frac{\partial}{\partial \bar z_m}(\overline{L}_{\nu}(\xi)f_{\tau}) \right\|_{\D_j} \\
&=&\left\| \frac{\partial}{\partial \bar z_m}(\overline{L}_{\nu}(\xi ))f_{\tau}\right\|_{\D_j}+\left\| \overline{L}_{\nu}(\xi )\frac{\partial}{\partial \bar z_m}(f_{\tau})\right\|_{\D_j}\\
&\leq& C(\|\dbar f\|_{\D_j}+\|\dbar^*f\|_{\D_j})<\infty
\end{eqnarray*}
for $m=1,2$ and a positive constant $C$ that does not depend on $f$. In the last inequality we used the fact that  $L^2$-norms of $f$ and the ``bar" derivatives of $f_{\tau}$ on $\D_j$ are bounded by $C(\|\dbar f\|_{\D_j}+\|\dbar^*f\|_{\D_j})$. We remark that it is essential that $\xi$ is complex valued and $\D$ is smooth in a neighborhood of $\overline{\D}_j$. Therefore, we showed that $\xi f\in Dom(\Box^j)$.

As for $(1-\xi)f$ being in $Dom(\Box)$. Since $\xi\equiv 1$ in a neighborhood of the boundary points of $\D_j$ that meets the band created using $\gamma_j$ and $\gamma_{j-1}$ we have $(1-\xi)f\equiv 0$ on $\D\setminus \D_j$. Also since $\overline{L}_{\nu}(1-\xi)=-\overline{L}_{\nu}(\xi)$ similar calculations as before show that  $(1-\xi)f\in Dom(\Box)$.
This completes the proof of the claim. 

\smallskip

We will use generalized constants in the sense that $\|A\|_{s,\D_j}\lesssim \|B\|_{s,\D_j}$ means that there is a constant $C=C(s,\D_j)>0$ that depends only on $s$ and $\D_j$ but not on $A$ or $B$  such that $\|A\|_{s,\D_j}\leq C\|B\|_{s,\D_j}$.
Assume that the $\dbar$-Neumann operator on $\D$ maps $W^s_{(0,1)}(\D)$ into itself continuously for some $s>0$. That is, $\|N_{1}h\|_{s,\D} \lesssim \| h\|_{s,\D}$ for $h\in W^s_{(0,1)}(\D)$. Then we have $\|g\|_{s,\D} \lesssim \| \Box g\|_{s,\D}$ for $g\in Dom(\Box)$ and $\Box g\in W^s_{(0,1)}(\D) $. 
Let $f\in Dom(\Box^j)$ and $\Box^j f\in W^s_{(0,1)}(\D_j) $. Then we have 
$$\|f\|_{s,\D_j}\leq \|\xi f\|_{s,\D_j}+\|(1-\xi )f\|_{s,\D_j}.$$
Since $\xi\equiv 0$ in a neighborhood of the weakly pseudoconvex boundary points of $\D_j$ we can use pseudolocal estimates on $\D_j$ (see \cite{KN}) to get
\begin{equation} \label{pseudolocal}
\|\xi f\|_{s,\D_j} \lesssim \|\Box^j f\|_{s-1,\D_j}+\|\Box^j f\|_{\D_j}.
\end{equation}
Let us choose $\eta$ to be a smooth compactly supported function that is constant $1$ around the support of $\nabla \xi$ and zero in a neighborhood of the weakly pseudoconvex points of $\D_j$. 
Therefore, we have 
\begin{eqnarray*}
\|(1-\xi)f\|_{s,\D_j} = \|(1-\xi)f\|_{s,\D} & \lesssim & \|\Box (1-\xi)f\|_{s,\D} \\
&\lesssim & \| (\bigtriangleup\xi )f\|_{s,\D}+\|\nabla\xi \cdot \nabla f\|_{s,\D}+\|(1-\xi)\bigtriangleup f\|_{s,\D_j}\\
&\lesssim & \|\eta f\|_{s,\D_j}+\|\eta f\|_{s+1,\D_j}+ \|\Box^j f\|_{s,\D_j} \\
&\lesssim & \|\Box^j f\|_{s,\D_j}.
\end{eqnarray*}
The first inequality comes from the assumption that the $\dbar$-Neumann operator on $\D$ is continuous on $W^s_{(0,1)}(\D)$. The second inequality comes from the fact that $\Box$ operates as Laplacian componentwise on forms. In the last inequality we used the pseudolocal estimates as we did in (\ref{pseudolocal}). 
Therefore we showed that 
$$\| f\|_{s,\D_j} \lesssim \|\xi f\|_{s,\D_j} + \|(1-\xi) f\|_{s,\D_j} \lesssim \|\Box^j f\|_{s,\D_j} $$ for $f\in Dom(\Box^j)$ and $\Box^jf\in W^s_{(0,1)}(\D_j)$. One can check that this is equivalent to the condition that the $\dbar$-Neumann operator on $\D_j$ is continuous on $W^s_{(0,1)}(\D_j)$.
\end{proof}

One can check that $\dbar^{*}N_{1}$ maps $W^s_{(0,1)}(\D)$ into $W^s(\D)$ continuously if and only if $\| \dbar^{*}f\|_{s,\D} \lesssim \|\Box f\|_{s,\D} $ for $f\in Dom(\Box)$ and $\Box f\in W^s_{(0,1)}(\D)$. 
Using this observation one can give a proof, similar to the proof of the theorem, for the following corollary. 

\begin{corollary}
There exists a bounded pseudoconvex domain $\D$ in $\C^2$, smooth except one point, such that 
$\dbar^*N_1$  is not bounded from $W^s_{(0,1)}(\D)$ into $W^s(\D)$ for any $s>0$. 
\end{corollary}

It is interesting that for a smooth bounded pseudoconvex domain $\D$ in $\C^2$ the operator $\dbar N_1$ is bounded from $W^s_{(0,1)}(\D)$ into $W^s_{(0,2)}(\D)$ for any $s\geq 0$. (One can use (4) in \cite{BS90} to see this). 

\begin{remark}
We would like to note the following additional property for  the domain we constructed in the proof of Theorem \ref{thm1}. There is no open set $U$ that contains the non-smooth boundary point of $\D$ such that $\overline{U \cap \D}$ has a Stein neighborhood basis. That is, non-smooth domains may not have a ``local" Stein neighborhood basis. However, this is not the case for smooth domains (see for example \cite[Lemma 2.13]{R86}). 

\end{remark}

\section*{Acknowledgment} I would like to thank Professors Harold Boas and Emil Straube for their encouragement and help for so many years, and Daniel Jupiter for teaching me the barbell lemma. I also would like to thank the referee for helpful suggestions.
\bibliographystyle{amsalpha}
\bibliography{regu}
\end{document}